\documentstyle[draft, amscd, syntonly, psfig, amssymb]{amsart}
\chardef\bslash=`\\
\newcommand{\ntt}{\series m\shape n\tt}

\newcommand{\cs}[1]{{\protect\ntt\bslash£1}}

\newcommand{\opt}[1]{{\protect\ntt£1}}

\newcommand{\env}[1]{{\protect\ntt£1}}
\def \s{\sigma}
\def \x{\times}
\def \sx{\s^{\x}}
\def \sp{\s^{+}}
\def \sn{\s^{-}}

\def \vu{\vec{u}}
\def \K{{\cal{K}}}
\def \B{{\cal{B}}}
\def \T{{\cal{T}}}
\def \L{{\cal{L}}}

\makeatletter
\def\verbatim{\interlinepenalty\@M \@verbatim
  \leftskip\@totalleftmargin\advance\leftskip2pc
  \frenchspacing\@vobeyspaces \@xverbatim}
\makeatother
\newtheorem{theo}{Theorem}[section]

\newtheorem{defi}{Definition}[section]
\newtheorem{remark}{Remark}[section]

\begin{document}

\title[Vassiliev Invariants: A formula]
{A formula in the theory of finite type invariants}

\author{N. A.  Askitas}
\address[N. Askitas]{. \newline Aristoteleion University of 
Thessaloniki, Mathematics Department, GR-54006, Thessaloniki, Greece.}
\email[N.Askitas]{askitas@@ccf.auth.gr and askitas@@mpim-bonn.mpg.de}

\date{}
\thanks{{\bf Acknowledgments}: Supported in part by EU TMR Fellowship 
No: ERBFMBICT982870, the author would like to thank Prof. G. Stamou 
for his help during his year at the Aristoteleion University of 
Thessaloniki, Greece and the Max-Planck-Institut f\"ur Mathematik
in Bonn for hospitality during the same time. He would also like to 
thank Hitoshi Murakami an exchange of emails with whom introduced him 
to Habiro's Theorem.}

 

\begin{abstract}  A family $C_{k+1}$ of local moves  on knot diagrams for 
each positive integer $k$,
is defined in \cite{r:kazuo} where it is also shown that two knots 
are $C_{k+1}$-equivalent iff all of
their Vassiliev-Gusarov invariants of degree $k$ agree. Every move 
$C_{k+1}$ splits the space $\K$ of knots into equivalence classes and 
defines a metric on each equivalence class. For two 
$C_{k+1}$-equivalent knots $K,J$
with $d_{C_{k+1}}(K,J)=1$  we give a formula for the difference 
$v_{k+1}(K)-v_{k+1}(J)$. From this we deduce a formula for the difference of 
the degree $k+1$ invariants of two knots all of whose degree $k$ 
Vassiliev invariants coincide.
\end{abstract}

\maketitle

\section{\bf Introduction}
We present a formula which expresses the difference of the degree 
$k+1$ Vassiliev invariants of two knots $K$ and $J$ which differ by 
one application of a certain local move\footnote{This paper grew out 
of part of a talk I gave in April of 1999 at the Topology Seminar of 
the Max-Planck-Institut f\"ur Mathematik in Bonn.}. In order to be 
able to state the result we need some terminology which we now 
explain.
\subsection{\bf Tangles, braids and symmetries}

Consider a quadrangle with $n$ points on its floor-side labeled 
$p_1,\ldots p_n$ and $n$
points on its ceiling-side labeled $p^1,\ldots ,  p^n$ as in the 
figure below:

{\vspace{.03in}}
\begin{figure}[h]
\centerline{\psfig{figure=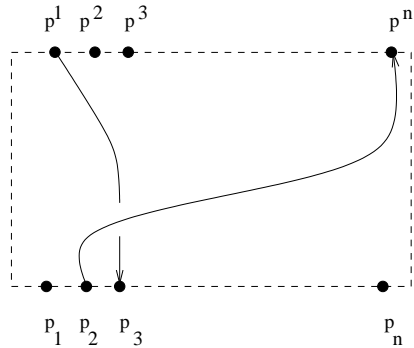,height=1.8in, clip=}}
\caption{The boundary points of braids in the $n$th braid group with 
an $o$-orientation and some connecting oriented strings.}
\end{figure}
{\vspace{.03in}}

These points will be the boundary points of oriented tangles after we 
introduce one
more piece of data.  Let $I_n = \{1, \ldots , n \}$ and let $o \colon 
I_n \rightarrow
Z_2$ be any function which orients the points $p_i ,p^j$ as follows: 
If $o(i)=0$ then
$p_i$ is an incoming point and $p^i$ is an outgoing point while if 
$o(i)=1$ then $p^i$ is an incoming point and $p_i$ is an outgoing 
point (in the figure above $o(2)=0=o(n)$ and $o(1)=o(3)=1$). 
Considering
oriented tangles whose strings connect a $p_i$ to a $p^j$, if 
$o(i)=0=o(j)$ or a $p^j$
to a
$p_i$, if $o(i)=1=o(j)$, gives rise to the $n$th $o$-braid group 
$\B_n^o$ whose end
points are oriented according to $o$. 

$\B_n^o$  can be thought of as follows. Let $\Phi_n: \B_n \rightarrow 
S_n$ be the well known group homomorphism which sends a braid to the 
permutation it induces on the boundary points. If we use $o$ to 
impose a partition on $I_n$ and then take the subgroup $S_n^o$ of 
$S_n$ which respects this partition then $B_n^o$ can be thought of as 
$\Phi_n^{-1}(S_n^o)$. 

Next consider general tangles with the same
boundary points as the ones in $\B_n^o$ (so as to for example allow 
connecting $p_i$
to $p_j$ when $o(i)+o(j)=1 \in Z_2$). Denote the set of all such 
tangles then by
$\T_n^o$ so as to have $\B_n^o \subset \T_n^o$. 

As a generalization of the group homomorphism $\Phi_n \colon \B_n 
\rightarrow S_n$ there is a map $\Phi_n^o
\colon \T_n^o \rightarrow S_n$ defined by sending a tangle $t \in 
\T_n^o$ to a
bijection $\Phi_n^o(t) \colon \{ p_i \colon o(i)=0 \} \cup \{ p^i 
\colon o(i)=1 \}
\rightarrow 
\{ p_i \colon o(i)=1 \} \cup \{ p^i \colon o(i)=o \}$ which sends $x 
\in \{ p_i \colon
o(i)=0 \} \cup \{ p^i \colon o(i)=1 \}$ to $y \in \{ p_i \colon 
o(i)=1 \} \cup \{ p^i
\colon o(i)=o \}$ if $t$ contains an oriented string connecting $x$ 
to $y$ (resp. $y$
to $x$) if $o(x)=0=o(y)$ (resp. $o(x)=1=o(y)$). Another way to say 
this is that oriented strings define a bijection from incoming to 
outgoing points. To pass from this bijection to an
element of $S_n$ we need to enumerate the points of the domain and 
range of
$\Phi_n^o$. We do this from left to right. The restriction then of 
this map $\Phi_n^o$
to $\B_n^o$ gives rise to a group homomorphism. 

By closure of elements of
$\T_n^o$   we will mean the identification of $p_i$ with $p^i$ for each 
$i$ as in the figure below. For $t \in \T_n^o$ its closure will be 
denoted by $\overline{t}$. The closure of an element  $t \in \T_n^o$ 
is a knot iff $\Phi_n^o(t)$ is an $n$-cycle. The operation of closure 
defines a map from $\T_n^o$ to $\L$ the space of link types. 

{\vspace{.03in}}
\begin{figure}[h]
\centerline{\psfig{figure=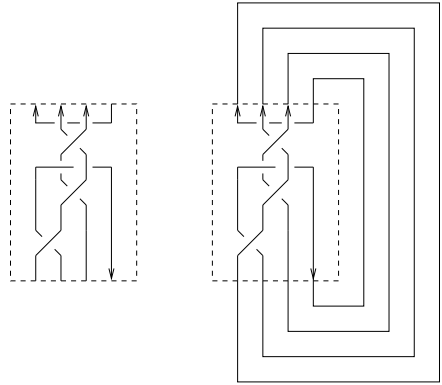,height=2in, clip=}}
\caption{Left: An $o$-braid in $\B_4^o$. Right: Its closure. Here 
$o\colon I_4 \rightarrow Z_2$ is given by: $o(1)=o(2)=o(3)=o(4)+1 =0 
\in Z_2$.}
\end{figure}
{\vspace{.03in}}


Now let $B_{n,e}^o={\Phi_n^o}^{-1}(e\in S^n)$ be the pure $n$th 
$o$-braid group. We
have a map: $B_{n,e}^o \times \T_n^o \rightarrow \T_n^o$ which sends 
$(b,t) \in
B_{n,e}^o \times \T_n^o$ to $bt \in \T_n^o$ by placing $t$ on top of 
$b$ as in the
figure below. This map satisfies the condition 
$\Phi_n^o(bt)=\Phi_n^o(t)$.

{\vspace{.03in}}
\begin{figure}[h]
\centerline{\psfig{figure=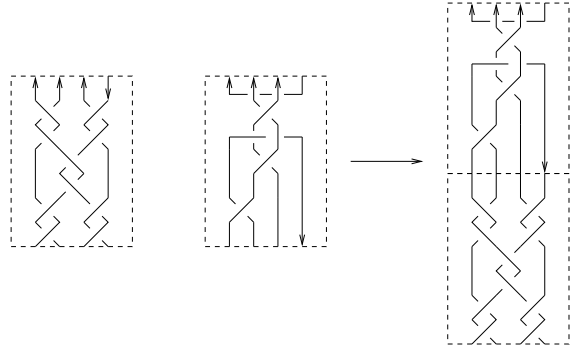,height=1.8in, clip=}}
\caption{ Left: A pure braid in $\B_4^o$. Center: A braid in 
$\T_4^o$. Right: Their composition.}
\end{figure}
{\vspace{.03in}}

\subsection{\bf Verginian moves}
In this section we introduce terminology which on the one hand 
motivates the main theorem and on the other makes its statement 
concise. We formalize to a certain extend the notion of local move on 
knot diagrams.

\begin{defi} \label{d:vergina}
A {\bf Verginian} move or operator\footnote{These are called local 
moves in the literature and
sometimes Gordian moves as R. Wendt coined them. However in the story 
involving
Alexander the Great, Gordian is the knot and not the cutting method. 
So one should
call such a move a Great Alexandrian move. Due however to the 
existence of a great
Alexander in topology I chose to refer to them as Verginian 
operations after the
birthplace of Alexander the Great.}  on the space $\K$ of all knots 
is a pair of
tangles $(T_1, T_2)$ with the same boundary acting on a knot $K$ by 
scanning a
projection of $K$ for the appearance of $T_1$ (resp. $T_2$) and then 
replacing it by
$T_2$ (resp. $T_1$). We say that $K$ and $J$ can be connected via 
$\mu$ moves if there
is a finite sequence of applications of $\mu$ starting with a 
projection of $K$ and
ending with one of $J$. In this manner we get an equivalence relation 
$\sim_{\mu}$ on $\K$ and
on each equivalence class we get a metric (we denote them all by 
$d_{\mu}$) defined as
follows: If $K \sim_{\mu} J$ define $d_{\mu}(K,J)$ to be the minimum 
number of times
we need to apply $\mu$ to pass from a projection of $K$ to one of 
$J$.  The number $|\K / \sim_{\mu}|-1$ is called the {\bf 
unknotting deficiency} of the move $\mu$. A Verginian move
is called an unknotting operation if $\K / \sim_{mu}= \{ [O] \}$  
i.e. if its unknotting deficiency is $0$. 
\end{defi}

For every Verginian operator we get numerical knot invariants on every
$\mu$-equivalence class by fixing one base knot for every element of 
$\K /
\sim_{\mu}$ and taking the distance of any other knot in the same 
$\mu$-class to that
base knot. The ordinary unknotting number uses as base knot the 
unknot.  Any numerical
invariant thus obtained is of non- finite type. The issue therefore 
arises of how (if
at all) such numerical invariants are encoded into the finite type 
ones especially in
view of the conjectured knot classification by finite type invariants.

Verginian moves can be composed as follows. Let $\mu_i$, $i=0,1$ be 
two such moves. Denote their constituting coordinate tangles by 
$p_j(\mu_i)$ where $p_j$ is the projection on the $j$th factor 
$j=1,2$. Then the union move $\mu_1 \cup \mu_2$ is that move which 
acts on knots by scanning knot projections for the appearance of 
$p_j(\mu_i)$ and replacing it by $p_{j+1} (\mu_i)$, where $j$ is to 
be read modulo $2$.

\subsection{\bf The Definition of the move and the statement of the 
Theorem}

Now we are ready to define the move we want to consider and state the 
main Theorem of the paper.
 Let
$C_{k,\vec{d},o}$ be the local move defined by $(BH_{\vec{d}}^o(k),e) 
\in \B_{k+2}^o
\times B_{k+2}^o$, where 

$$BH_{\vec{d}}^o(k)= (\prod_{i=1}^{k}
\s_i^{d_i})\s_{k+1}^{d_{k+1}}(\prod_{i=1}^{k}\s_{k+1-i}^{-d_{k+1-i}})
(\prod_{i=2}^{k}\s_{i}^{d_i})\s_{k+1}^{-d_{k+1}}
(\prod_{i=2}^{k}\s_{k+2-i}^{-d_{k+2-i}}),$$

with  $d_{i} \in \{ \pm 2\}$, $\vec{d}=(d_1, \ldots , d_{k+1})$ and 
$\s_i$ are the standard generators of $B_{k+2}^o$ (here $o \colon 
I_{k+2} \rightarrow Z_2$). This is an element of $B_{k+2}^o$ and we 
will use the same symbol to denote the geometric tangle defined by it.

\begin{remark} \label{rem:habiro}

The reader may check that the standard closure of $BH_{\vec{d}}^o(k)$ 
in nothing but the $k$th iterated Bing-double of the Hopf-link, for 
any $\vec{d}, o$ and that therefore the union
$\cup_{o, \vec{d}}C_{k,o, \vec{d}}$ is nothing but K. Habiro's 
$C_{k+1}$-move (\cite{r:kazuo}). In the language of Verginian moves: 
$$C_{k+1}=\bigcup_{o, \vec{d}}C_{k,o, \vec{d}}$$

\end{remark}

{\vspace{.03in}}
\begin{figure}[h] \label{f:habiro}
\centerline{\psfig{figure=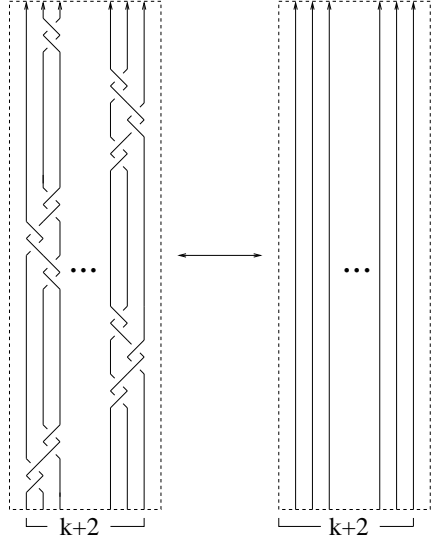,height=2.8in, clip=}}
\caption{The move $C_{k,\vec{d},o}$ for $d_i=2$ for all $i$ and ${\rm 
Im}(o)={0}$ }
\end{figure}
{\vspace{.03in}}
We will write $d_{{k,\vec{d},o}}$ instead of $d_{C_{k,\vec{d},o}}$ 
for the metrics $C_{k,\vec{d},o}$ defines on each of its equivalence 
classes.

Now we are ready to state the theorem of this paper after we define 
some braids in $B_{k+2}^o$. Let $\vu \in Z_2^{k}$ be a vector 
parameter and denote its coordinates by $u_i$. For any such $\vu$ 
define $\vu +1$ to be the vector $ \vu +(1,1, \ldots , 1)$.  Define 
then $W_{\vu}=a_1 \cdots a_k$ with $ a_i=e$ or $\s_i^2$ according as 
$u_i$ is $0$ or $1$.  Define also $W_{\vu}^{r}=a_k \cdots a_1$ with $ 
a_i=e$ or $\s_i^2$ according as $u_i$ is $0$ or $1$. 

\begin{theo} \label{th:main}
Suppose that $K$ and $J$ are two knots with 
$d_{{k,\vec{d},o}}(K,J)=1$ for some $\vec{d}, o, k$. Then we can write $K$ 
as the closure $\overline{BH_{\vec{d}}^o(k)T}$ of $BH_{\vec{d}}^o(k)T$ and 
$J$ as the closure $\overline{T}$ of $T$ for some $T \in \B_{k+2}^o$. 
Furthermore: 

$$ v_{k+1}(K)-v_{k+1}(J) = s(d_{k+1})(-1)^{o(k+1)o(k+2)}\sum_{\vu \in 
Z_2^k} s_{\vu , \vec{d}, o}
v_{k+1}(\overline{W_{\vu} \s_{k+1}^2 W_{\vu +1 }^{r}x})$$

for any Vassiliev invariant $v_{k+1}$ of degree $k+1$ and any $x \in \T_{k+2}^o$ such that 
$\Phi_{k+2}^o(T)=\Phi_{k+2}^o(x)$. The sign
$s_{\vu , \vec{d}, o}$ is given by:

$$ s_{\vu , \vec{d}, o}=\prod_{i=1}^{k} (-1)^{u_i+1}  
s(d_i)(-1)^{o(i)o(i+1)} ) $$

\end{theo}

\begin{remark} \label{rem:inspire}
For $k=1$ this is essentially the result in \cite{r:okada} from which 
this paper was originally inspired. During its conception the author 
was made aware of the newly published \cite{r:ohyama} where a result 
similar in philosophy is obtained. The result of this paper is 
different than that of \cite{r:ohyama} in terms of approach, degree 
of analysis and scope.

\begin{itemize}

\item Notice that because $d_{C_{k+1}}(K,J)=1$ implies 
$d_{{k,\vec{d},o}}(K,J)=1$ for some $\vec{d}$ and some $o$ the 
theorem states that if $d_{C_{k+1}}(K,J)=1$ then the formula of 
Theorem \ref{th:main} holds for some $\vec{d}$ and $o$ and any $x \in 
\T_{k+2}^o$ such that $\Phi_{k+2}^o(x)=\Phi_{k+2}^o(T)$. 

\item Notice that it follows from the formula of Theorem 
\ref{th:main} that, for fixed $K$, $J$, $o$, $\vec{d}$ with 
$d_{{k,\vec{d},o}}(K,J)=1$, the degree $k+1$ invariants do not see 
the cycle $\Phi_{k+2}^o(x)$. In fact it is messy but not hard to convince 
oneself that the knots which appear on the right hand side of the 
formula of Theorem \ref{th:main} when taken as a set with 
multiplicities remains invariant under the choice of $x$ in the 
following sense. Let $\s$ and $\s^{'}$ be two $k+2$-cycles in 
$S_{k+2}$. Then there exist tangles $x$ and $x^{'}$ with 
$\Phi_n^o(x)=\s$ and $\Phi_n^o(x^{'})=\s^{'}$ such that set of knots  
$\overline{ W_{\vu} \s_{k+1}^2 W_{\vu +1 }x}$ is equal to the set of 
knots  $\overline{W_{\vu} \s_{k+1}^2 W_{\vu +1 }x^{'}}$.

\end{itemize}

\end{remark}

{\vspace{.03in}}
\begin{figure}[h]
\centerline{\psfig{figure=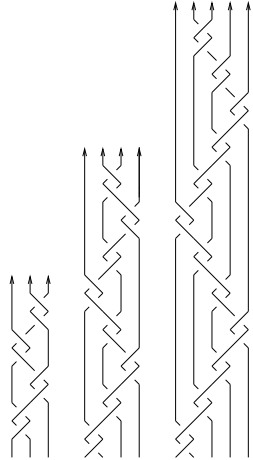,height=2.4in, clip=}}
\caption{$BH_{\vec{d}}^o(k)$ for $k=1,2,3$, ${\rm{Im}}(o)=0$, $d_i=2$ 
all $i$.  }
\end{figure}
{\vspace{.03in}}

If two knot $K$ and $J$ differ by replacing one of the braids above 
with the corresponding trivial one then the difference of their 
degree $k$ (=2,3,4) invariants is given as in the figure below. We 
assume for the sake of simplicity and concreteness that $T$ maps via 
$\Phi_{k+2}$ to the $k+2$ cycle $(1 (k+2) (k+1) k \ldots  2 )$.

{\vspace{.03in}}
\begin{figure}[h]
\centerline{\psfig{figure=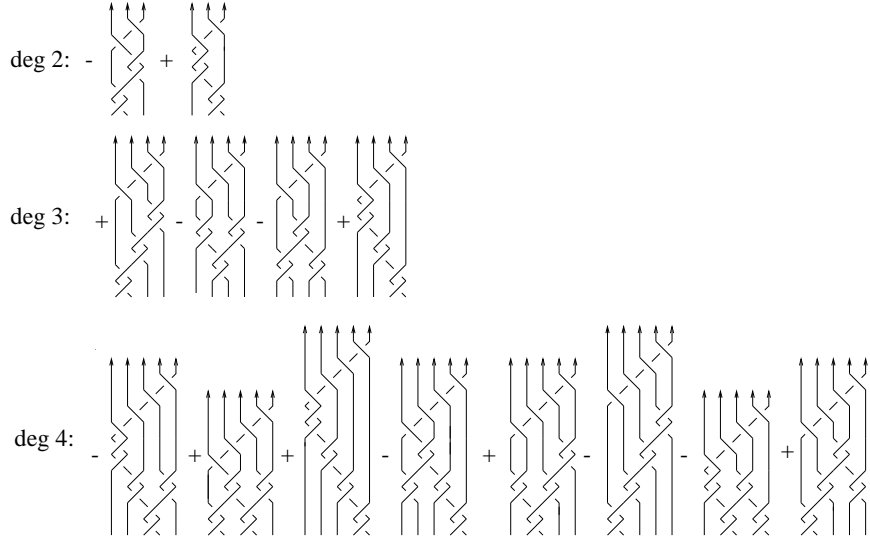,height=2.8in, clip=}}
\caption{ The difference of Vassiliev invariants. For every tangle 
get a knot by closing. In place of every knot the value of its 
respective invariant. }
\end{figure}
{\vspace{.03in}}

\section{\bf The proof}
We introduce some notation and remind the reader of standard 
conventions. Once we are set with all of that the proof will 
be straightforward. 
We adopt the convention that when we write equations of Vassiliev 
invariants with several (singular) knots involved we draw (or write 
whichever the case maybe) only that part of the knot where the knots 
may differ. Also in an equation of Vassiliev invariants we can write 
a knot $K$ instead of the value $v(K)$ of the invariant $v$ on $K$ if 
$v$ is understood. We will be writing the generators of $B_n^o$ as 
$\s_i$. There is an ambiguity of notation here but we will avoid it 
by declaring $o$ every  time we write $\s_i$'s. We will also consider 
singular $o$-braids. We will write $\sx_i$ for the case where we have 
a self intersection at the corresponding place. In the figure below 
we can see an example of this notation.

{\vspace{.03in}}
\begin{figure}[h]
\centerline{\psfig{figure=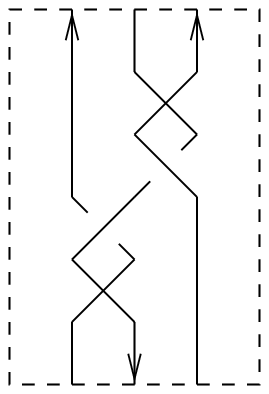,height=.7in, clip=}}
\caption{ A $2$-singular $o$-braid on $3$ strings. Here $o\colon I_3 
\rightarrow Z_2$ is given by $o(1)=o(3)=o(2)+1=0 \in Z_2$. This 
singular braid can be written as: $\sx_1 \s_1 \s_2^{-1} \sx_2$.}
\end{figure}
{\vspace{.03in}}

We will use the notation $\sp_i$ instead of $\s_i \sx_i = \sx_i \s_i$ 
and $\sn_i$ instead of $\s_i^{-1} \sx_i = \sx_i \s_i^{-1}$ so that 
the example in the figure above can be written as $\sp_1 \sn_2$. For 
any non-zero real number $x$ we will write $s(x)$ for 
$\frac{|x|}{x}$. Define $s_x=\pm$ whenever $s(x)= \pm 1$. Instead of 
$s_{d_i}$ we will write $s_i$. We need this symbol in order to be 
able to express ambiguity of the type $\s_i^{\pm}= \sx_i \s_i^{\pm 
1}$. 
We draw the readers attention to the difference between $\s_i^{\pm}$ 
and $\s_i^{\pm 1}$. The first are singular words the second ones 
non-singular.
Then if a knot projection contains $\s_i^{d_i} \in \B_n^o$ 
we can express the Birman-Lin condition (\cite{r:birlin}) at a crossing of $\s_i^{d_i}$ by writing: $ \s_i^{d_i} = e + s(d_i) (-1)^{o(i)+o(i+1)} \s^{s_i}$ or simply  $ \s_i^{d_i} = e + so(i) \s^{s_i}$ if we let $so(i)=s(d_i) 
(-1)^{o(i)+o(i+1)}$.

Now suppose that we are given a word in some $B_n^o$ and we want to 
apply the Birman-Lin condition:

{\vspace{.03in}}
\begin{figure}[h]
\centerline{\psfig{figure=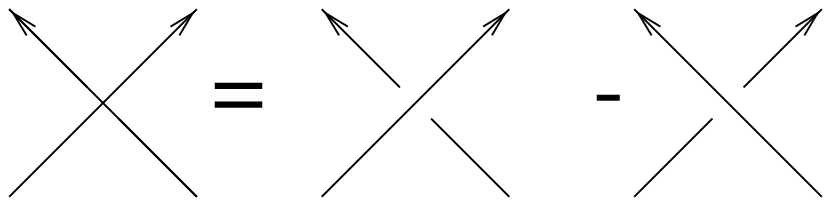,height=.3in, clip=}}
\caption{ The Birman-Lin Condition.}
\end{figure}
{\vspace{.03in}}

Let's say that we have $\s_1^2 \s_2^2 \s_1^{-2} \s_2^{-2} \in \B_3$ and 
that we have marked a crossing where we will apply the Birman-Lin 
condition. We can then write out the tree of possibilities and at the 
end look at our expression or we can write the tree in the form of an 
expression like: $(e+\sp_1) (e+\sp_2) (e-\sp_1) (e-\sp_2)$ and then just 
multiply through. The reader may check for himself that such a product 
is meaningful as well as check the steps we give below.

\begin{equation}
\begin{array}{l}
(e+\sp_1) (e+\sp_2) (e-\sn_1) (e-\sn_2) =\\
(e+\sp_1) e (e-\sn_1) (e-\sn_2) +(e+\sp_1) \sp_2 (e-\sn_1) 
(e-\sn_2)=\\
(e+\sp_1)(e-\sn_1)(e-\sn_2)+(e+\sp_1) \sp_2 (e-\sn_1) (e-\sn_2)= \\
(e+\sp_1-\sn_1 -\sp_1 \sn_1) (e-\sn_2)+(e+\sp_1) \sp_2 (e-\sn_1) 
(e-\sn_2)=
\end{array}
\end{equation}

but now notice that $\sp_1 \sn_1 = \sx_1 \sx_1= \sp_1 - \sn_1$ so we 
can continue as follows:

\begin{equation}
\begin{array}{l}
(e-\sn_2)+(e+\sp_1) \sp_2 (e-\sn_1) (e-\sn_2)= \\
e-\sn_2+ \sp_2  (e-\sn_2)
+\sp_1 \sp_2  (e-\sn_2)= \\
-\sp_2 \sn_1 (e-\sn_2)
-\sp_1 \sp_2 \sn_1 (e-\sn_2) =\\
e-\sn_2+ \sp_2 -\sp_2  \sn_2
+\sp_1 \sp_2  -\sp_1 \sp_2  \sn_2= \\
-\sp_2 \sn_1 +\sp_2 \sn_1 \sn_2
-\sp_1 \sp_2 \sn_1 + \sp_1 \sp_2 \sn_1 \sn_2 =\\
e+\sp_1 \sp_2-\sp_2 \sn_1 +O_3 
\end{array}
\end{equation}

where $O_3$ contains singular braids with at least three 
singularities.
So if this was a calculation of $v_2$ then $O_3$ would vanish and 
hence we would have Okada's result in \cite{r:okada}. Examining this 
calculation the reader may get familiar with this notation quicker 
than it would take to write out all the details. The idea is that if 
we have a knot projection containing a braid  and we intend to apply 
the Birman-Lin condition only locally on the braid in order to 
express its Vassiliev invariant of some degree in terms of singular 
knots resulting from these considerations and we therefore adopt the 
convention that we only write the braid instead of the invariant then 
we can perform the calculations by treating singular braids as 
objects in an algebra of sorts.

 We now proceed to proving Theorem \ref{th:main} in its complete 
generality.

Suppose that $K$ and $J$ are two knots such that 
$d_{{k,\vec{d},o}}(K,J)=1$ for some $\vec{d},o$. Then we can write 
$K$ as $\overline{BH_{\vec{d}}^o(k)T}$ the closure of 
$BH_{\vec{d}}^o(k)T$ and $J$ as $\overline{T}$ for some $T \in 
\T_{k+2}^o$.

We now compute the degree $k+1$ Vassiliev invariant of $K$. Because 
we will only use the Birman-Lin Condition on $BH_{\vec{d}}^o(k)$ we 
can write $v_{k+1}(BH_{\vec{d}}^o(k))$ instead of $v_{k+1}(K)$ and in 
fact since we have declared that we are computing $v_{k+1}$ we can 
drop the $v_{k+1}$ all together and write $BH_{\vec{d}}^o(k)$ instead 
of $v_{k+1}(BH_{\vec{d}}^o(k))$.

\begin{equation}
\begin{array}{l}
BH_{\vec{d}}^o(k) =\\
\displaystyle{\prod_{i=1}^{k}}(e+so(i)\s_i^{s_i}) 
(e+so(k+1)\s_{k+1}^{s_{k+1}}) 
\displaystyle{\prod_{i=1}^{k}}(e-so(k+1-i)\s_{k+1-i}^{-s_{k+1-i}})\\
\displaystyle{\prod_{i=2}^{k}}(e+so(i)\s_{i}^{s_{i}})(e-so(k+1)\s_{k+1}^{-s_{k+1}})
\displaystyle{\prod_{i=2}^{k}}(e-so(k+2-i)\s_{k+2-i}^{-s_{k+2-i}})
\end{array}
\end{equation}
We will expand this expression in two ways. Once starting from 
$(e+so(k+1)\s_{k+1}^{s_{k+1}})$ and proceeding radially in both 
directions and once from $(e\pm so(1)\s_1^{\pm s_1})$ Comparing the 
two expressions will give us the desired result. The key relation is 
this: $\sp_i \sn_i =\sx_i \sx_i =
\sp_i -\sn_i$ by means of which we have:
\begin{equation} \label{eq:key}
(e+so(i)\s_i^{s_i}) (e-so(i)\s_i^{-s_i})= e + so(i)\s_i^{s_i} 
-so(i)\s_i^{-s_i} -\s_i^{s_i}\s_i^{-s_i}=e.
\end{equation}
We begin by expanding radially from  the peak 
$(e+so(k+1)\s_{k+1}^{s_{k+1}})$ down the hills on both sides.
Before we begin we outline the plan. We will expand with an eye on 
separating the final expression into sums of singular words which 
contain $\s_1^{\pm }$ and ones which do not. We will now argue that 
$BH_{\vec{d}}^o(k)$ is given by:

\begin{equation}\label{eq:once}
BH_{\vec{d}}^o(k)= e+ \sum_{\vu \in Z_2^{k}} so_{\vu} U_{\vu} 
\s_{k+1}^{s_{k+1}}U_{\vu+1}^{-1}+ O
\end{equation}

where $W_{\vu}= b_1 \ldots b_n$  with $b_i= e$ or $\s_i^{s_i}$ 
depending on whether the $i$th coordinate of $\vu$ is $0$ or $1$. 
$W_{\vu}^{-1}$ is the inverse of the singular word  $W_{\vu}$ (upside 
down crossing reversed singularities follow accordingly). $O$ is a 
sum of singular words which do not contain $\s_1^{\pm}$. The sign 
$so_{\vu}$ is given by:

$$ so_{\vu}= \prod_{i=1}^{k} so(i) (-1)^{u_i+1}$$ 

To see that equation (\ref{eq:once}) holds is not hard. It falls 
basically out of the relations of equation (\ref{eq:key}) as follows. 

Let $$A= \displaystyle{\prod_{i=1}^{k}}(e+so(i)\s_i^{s_i}) 
(e+so(k+1)\s_{k+1}^{s_{k+1}}) 
\displaystyle{\prod_{i=1}^{k}}(e-so(k+1-i)\s_{k+1-i}^{-s_{k+1-i}}),$$
and begin expanding it from the middle outward. Then using equation 
(\ref{eq:key}) we get:

$$A= e+ so(k+1)\displaystyle{\prod_{i=1}^{k}}(e+so(i)\s_i^{s_i}) 
\s_{k+1}^{s_{k+1}}
\displaystyle{\prod_{i=1}^{k}}(e-so(k+1-i)\s_{k+1-i}^{-s_{k+1-i}})$$

Now notice that:

\begin{equation} \label{eq:step1}
\begin{array}{l}
so(k+1)\displaystyle{\prod_{i=1}^{k}}(e+so(i)\s_i^{s_i}) 
\s_{k+1}^{s_{k+1}}
\displaystyle{\prod_{i=1}^{k}}(e-so(k+1-i)\s_{k+1-i}^{-s_{k+1-i}})= \\
 \\
so(k+1)\s_{k+1}^{s_{k+1}}+\\
so(k)so(k+1)\displaystyle{\prod_{i=1}^{k-1}}(e+so(i)\s_i^{s_i}) 
\s_{k}^{s_{k}} \s_{k+1}^{s_{k+1}}
\displaystyle{\prod_{i=2}^{k}}(e-so(k+1-i)\s_{k+1-i}^{-s_{k+1-i}}\\
-  so(k)so(k+1)\displaystyle{\prod_{i=1}^{k-1}}(e+so(i)\s_i^{s_i}) 
\s_{k+1}^{s_{k+1}}\s_{k}^{-s_{k}}
\displaystyle{\prod_{i=2}^{k}}(e-so(k+1-i)\s_{k+1-i}^{-s_{k+1-i}})\\
-so(k+1)\displaystyle{\prod_{i=1}^{k-1}}(e+so(i)\s_i^{s_i}) 
\s_{k}^{s_{k}}
 \s_{k+1}^{s_{k+1}} \s_{k}^{-s_{k}}
\displaystyle{\prod_{i=2}^{k}}(e-so(k+1-i)\s_{k+1-i}^{-s_{k+1-i}})
\end{array}
\end{equation}

Inspecting equation (\ref{eq:step1}) we can get equation 
(\ref{eq:once}) as follows. The middle two summands are the only ones 
which have a chance of leading to singular words containing 
$\s_1^{\pm}$. That the first summand does not it is clear because it 
does not contain $\s_1^{\pm}$ and cannot be expanded further. The 
last summand will not  lead to singular words containing $\s_1^{\pm}$ 
because the only way it has a chance to do this is through words 
where for each $i$ at least one of $\s_i^{\pm s_i}$ remains left or 
right of $\s_{k+1}^{s_{k+1}}$. This is because as soon as for some 
$i$ this does not happen the $\s_1^{\pm}$'s will disappear using 
equation (\ref{eq:key}).
But if this is the case since we are computing invariants of degree 
$k+1$ the terms which retain the chance of containing $\s_i^{\pm s_i}$ 
will disappear before the expansion reaches $\s_1^{\pm}$.

Now we begin expanding the expression for $BH_{\vec{d}}^o(k)$ in the 
second way.
Write  $BH_{\vec{d}}^o(k) = (e+so(1)\s_1^{s_1})M (e-so(1)\s_1^{-s_1}) 
N$ and expand to get:

 \begin{equation} \label{eq:twice}
\begin{array}{l}
BH_{\vec{d}}^o(k)= 
MN +so(1)\s_1^{s_1}MN -so(1)M \s_1^{-s_1} N - \s_1^{s_1}M \s_1^{-s_1} 
N
\end{array}
\end{equation}
Notice that using equation (\ref{eq:key}) $MN=e$. Also all the other 
terms contain $\s_1^{\pm}$. Comparing equations (\ref{eq:once}) and 
(\ref{eq:twice}) we conclude that $O$ is zero and that

\begin{equation}\label{eq:singular}
BH_{\vec{d}}^o(k)= e+ \sum_{\vu \in Z_2^{k}} so_{\vu} U_{\vu} 
\s_{k+1}^{s_{k+1}}U_{\vu+1}^{-1}.
\end{equation}

But now remember that (\ref{eq:singular}) is an equation of degree 
$k+1$ Vassiliev invariants of singular knots which are identical away 
from the singular words appearing in the equation. On the other hand 
in the equation every singular knot has $k+1$ singularities and hence 
it's crossings can be changed at will. This means that we can for 
example for each of this singular braids make any choice of a tangle 
$D_{\vu}$ with which to close them as long as $\Phi_n^o 
(D_{\vu})=\Phi_n^o(B)$. Also we can change the crossings of the 
singular braids so that instead of $U_{\vu+1}^{-1}$ we have 
$U_{\vu+1}^r= a_n \ldots a_1$ with $ a_i=e$ or $\s_i^{+}$ 
according as $u_i=1$ or $0$ respectively. Also we can assume that 
$U_{\vu}$ is given by $ a_1 \ldots a_n$ with $ a_i=e$ or $\s_i^{+}$ 
according as $u_i=0$ or $1$ respectively  So we can write equation 
(\ref{eq:singular}) as:

\begin{equation}\label{eq:sing}
BH_{\vec{d}}^o(k)= e+ \sum_{\vu \in Z_2^{k}} so_{\vu} 
\overline{U_{\vu} \s_{k+1}^{s_{k+1}}U_{\vu+1}^{r}D_{\vu}},
\end{equation}
for any choice of $D_{\vu}$ as above. If the choice is made so that 
$D_{\vu}=x$=const. for all $\vu$ then we can desingularize via the 
Birman-Lin Condition and find that all but the desired terms are 
canceled because they pair up with opposite signs on the level of 
braids already. \qed

\section{\bf Closing Remarks}

We would like now to see what Theorem \ref{th:main} says about the 
Vassiliev invariants of degree $k+1$ of $C_{k+1}$-equivalent knots 
i.e. (via Habiro's Theorem)
of knots all of whose degree $k$ invariants agree. Consider the 
$o$-braid group $\B_{n(k+2)}^o$ for $o \colon I_{n(k+2)} \rightarrow 
Z_2$ an orientation function, the $o$-tangle set $\T_{n(k+2)}^o$ and 
the function $\Phi_{n(k+2)}^o \colon \T_{n(k+2)}^o \rightarrow 
S_{n(k+2)}$. Consider also the pure braid group $\B_{n(k+2),e}^o= 
{\Phi_{n(k+2)}^o}^{-1}(e \in S_{n(k+2)})$.

In $\B_{n(k+2),e}^o= {\Phi_{n(k+2)}^o}^{-1}(e \in S_{n(k+2)})$ 
consider the words as in the figure below:
{\vspace{.03in}}
\begin{figure}[h]
\centerline{\psfig{figure=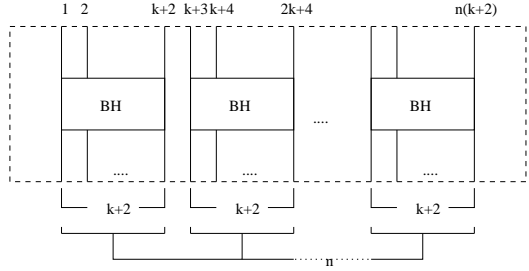,height=1.4in, clip=}}
\caption{ The words $BH_{\vec{d}}^o(n,k) \in \B_{n(k+2),e}^o$ formed 
by sticking sideways the words $BH_{\vec{d_j}}^{o_j}(k)$ (inside each 
box labeled $BH$) defined on the points $I_{j,k+2}=\{ j+1, \ldots 
j+k+2 \}$ for $j=0, \ldots , n-1 \}$. The orientation functions $o_j$ 
are given by restricting $o$ on $I_{j,k+2}$. }
\end{figure}
{\vspace{.03in}}

Form also words $W_{\vu_j} \in \B_{n(k+2),e}^o$ as in the figure 
below:

{\vspace{.03in}}
\begin{figure}[h]
\centerline{\psfig{figure=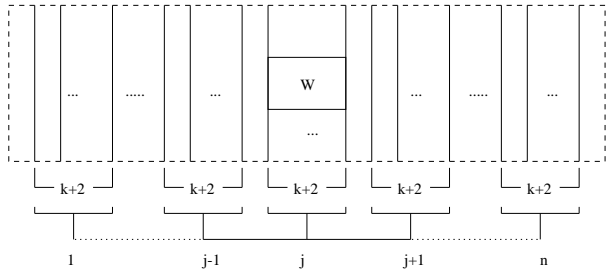,height=1.4in, clip=}}
\caption{ The words $W_{\vu_j} \in \B_{n(k+2),e}^o$ formed by means 
of the word $W_{\vu}$ (inside the box labeled $W$) defined on the 
points $I_{j,k+2}=\{ j+1, \ldots j+k+2 \}$ for $j=0, \ldots , n-1 
\}$. }
\end{figure}
{\vspace{.03in}}
Define now the (Habiro)-set of pure $o$-braids $\B_{n(k+2),H}^o$ 
containing all braids of the form $BH_{\vec{d}}^o(n,k)$. This is 
nothing but the preimage under the closure map (let us call it 
$\kappa$ here), $\kappa \colon \T_{n(k+2)}^o \rightarrow \L$, where 
$\L$ is the space of all links, of the $k$th iterated Bing-double of 
the Hopf link. Then Habiro's Theorem that ``two knots are 
$C_{k+1}$-equivalent iff all their degree $k$ Vassiliev invariants 
coincide" can be restated as saying that two knots are 
$C_{k+1}$-equivalent iff there is a $T_{n_k(k+2)} \in \T_{n_k(k+2)}^o$ 
with $\Phi_{n_k(k+2)}^o(T_{n_k(k+2)})$ a $n_k(k+2)$-cycle $ \in S_{n_k(k+2)}$ 
for some $n_k$ and $o$ such that while $J= \overline{T_{n_k(k+2)}}$, $K$ 
is in the image of the composite map below (where the first map is 
the composition defined in Figure 3)

$$\B_{n_k(k+2)}^o \times \{ T_{n_k(k+2)} \} \rightarrow \T_{n_k(k+2)}^o 
\stackrel{\kappa}{\rightarrow} \K$$

Now let 
$$ \B_{n_k(k+2),V}^o = \bigcup_{j, \vu_j} \{ W_{\vu_j} 
\s_{j(k+2)+k+1}W_{\vu_j}^r \} $$
be the (Vassiliev)-set of braids.
Consider the map:

$$ \B_{n_k(k+2),V}^o \times {\Phi_{n_k(k+2)}^o}^{-1}(T_{n_k(k+2)}) 
\rightarrow \T_{n_k(k+2)}^o \stackrel{\kappa}{\rightarrow} \K$$
and consider the sets of knots $\K_x$ obtained via this map for any 
choice of $x \in {\Phi_{n_k(k+2)}^o}^{-1}(T_{n_k(k+2)})$.

We can now say that if  $K$ and $J$ are $C_{k+1}$-equivalent then:

$$v_{k+1}(K)-v_{k+1}(J)= \sum_{K \in \K_x} s(K) v_{k+1}(K),$$
for any choice of $x$, where $s(K)$ is a sign; or we can write 
precisely:

\begin{equation} \label{eq:general}
v_{k+1}(K)-v_{k+1}(J)=\sum_{j=1}^{n_k-1} \sum_{\vu_j \in Z_2^k} 
so_{\vu_j}v_{k+1}(\overline{W_{\vu_j} \s_{j(k+2)+k+1}^2 W_{\vu_j}^r 
x}),
\end{equation}
and this for any $x \in  \T_{n_k(k+2)}^o$ which satisfies 
$\Phi_{n_k(k+2)}^o(x)=\Phi_{n_k(k+2)}^o(T_{n_k})$, where 
$W_{\vu_j}=a_{j(k+2)+1} \ldots a_{j(k+2)+k}$ with $a_{j(k+2)+1}=e$ or 
$\s_{j(k+2)+1}^{2} $ according as the $i$th coordinate $u_{ji}$ of 
$\vu_j$ is equal to $0$ or $1$ and $W_{\vu_j}^r=a_{j(k+2)+k} \ldots 
a_{j(k+2)+1}$. The $d_{ji}$'s take values in $\{ \pm 2 \}$. The signs 
$so_{\vu_j}$ are given by:

$$ so_{\vu}= \prod_{i=1}^{k} so(i) (-1)^{u_{ji}+1}$$

Here is an example of the general statement. Suppose that the $C_2$ 
distance of two knots $K$ and $J$ is equal to $2$. Then they may 
differ by a braid which looks like:

{\vspace{.03in}}
\begin{figure}[h]
\centerline{\psfig{figure=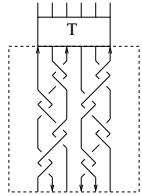,height=1in, clip=}}
\caption{ Two knots whose $C_2$ distance is equal to 2 may differ by 
such a tangle or its variants with respect to choices of $d_i$'s and 
and $o$'s (string orientations).}
\end{figure}
{\vspace{.03in}}

Then their degree $2$ Vassiliev Invariants differ by:
{\vspace{.03in}}
\begin{figure}[h]
\centerline{\psfig{figure=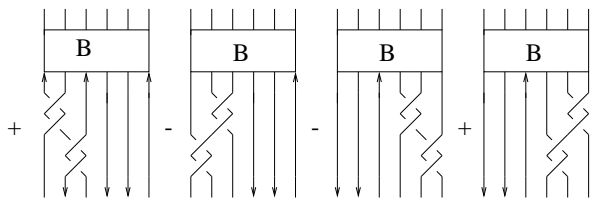,height=1in, clip=}}
\caption{ The difference of degree $2$ Vassiliev invariants. $B$ is 
any $o$-braid which induces the same full cycle in $S_6$ as $T$. Get 
knots by closing.}
\end{figure}
{\vspace{.03in}

A couple of final remarks:

\begin{itemize}
\item One of the consequences of Theorem \ref{th:main} is that the 
difference of the normalized degree $k+1$ Vassiliev invariants of 
knots with $d_{C_{k+1}}$-distance equal to one is bounded. More 
generally if we take the $d_{C_{k+1}}$-ball of radius $r$ then the 
difference of the normalized degree $k+1$ Vassiliev invariants of any 
two knots inside this ball is bounded by a constant which depends 
only on $r$.
\item It appears as though the permutation defined  by $T \in 
\T_{k+2}^o$ (the image under $\Phi_{k+2}^o$) may define a new 
invariant of knots which is not predicted by Vassiliev invariants. 
\item It seems as though the cardinality of the set of knots which appear in the formula of equation (\ref{eq:general}) should provide a (probably sharp) bound on the number of lineraly independent Vassiliev invariants of degree $k+1$. There appear to be many duplicates in their braid representation which when moded out should lead to such a bound (see Remark \ref{rem:inspire}).
We should also mention that it follows from the proof of Theorem \ref{th:main}
that the formula of the Theorem and its generalization hold also iff in the statement we replace all $\s_i^2$'s by $\s_i^{+}$'s. 
\end{itemize}

\end{document}